\documentclass{amsart}
\usepackage{amsmath}
\usepackage{amssymb}
\usepackage{amsthm}
\usepackage{amsfonts}
\usepackage{amstext}
\usepackage{amsopn}
\usepackage{amsxtra}
\usepackage{graphicx}
\usepackage{mathrsfs}

\newcommand{\ud}{\,{\rm d}}

\newcommand{\dt}{\frac{{\rm d}}{{\rm d}t}}

\newcommand{\D}{\mathrm{d}}

\newcommand{\T}{\mathbb{T}}

\newcommand{\be}{\begin{equation}}
\newcommand{\ee}{\end{equation}}

\DeclareMathOperator{\diverg}{div}

\newtheorem{lemma}{Lemma}[section]
\newtheorem{theorem}[lemma]{Theorem}
\newtheorem{proposition}[lemma]{Proposition}

\newtheorem{corollary}[lemma]{Corollary}
\theoremstyle{definition}

\theoremstyle{definition}
\newtheorem{remark}[lemma]{Remark}
\theoremstyle{definition}

{\catcode `\@=11 \global\let\AddToReset=\@addtoreset}
\AddToReset{equation}{section}

\newcommand{\N}{{\mathbb N }}
\newcommand{\R}{{\mathbb R}}

\newcommand{\ie}{{\sl i.e.\/ }}
\newcommand{\cf}{{\sl cf.\/ }}
\newcommand{\eg}{{\sl e.g.\/}}
\newcommand{\norm}[1]{ \left| \! \left| #1 \right| \! \right| }

\def\({\left(}
\def\){\right)}
\def\<{\left\langle}
\def\>{\right\rangle}

\renewcommand{\L}{{L}}
\renewcommand{\hbar}{{\displaystyle\bar{\phantom{x}}\kern-6pt h}}

\numberwithin{equation}{section}

\begin{document}


\title[Fokker-Planck equation for Bosons and Fermions]{Stability of steady states in kinetic 
Fokker-Planck equations for Bosons and Fermions} 
\author[L. Neumann]{Lukas Neumann}
\address{Johann Radon Institute for Computational and Applied Mathematics, Altenbergerstra\ss e 69, A-4040 Linz, Austria} 
\email{lukas.neumann@ricam.oeaw.ac.at}

\author[C. Sparber]{Christof Sparber}
\address{Wolfgang Pauli Institute Vienna \& Faculty of Mathematics, Vienna University, Nordbergstra\ss e 15, A-1090 Vienna, Austria}
\email{christof.sparber@univie.ac.at}

\begin{abstract}
We study a class of nonlinear kinetic Fokker-Planck type equations modeling quantum particles which obey the Bose-Einstein and 
Fermi-Dirac statistics, respectively. We establish the existence of classical solutions in the perturbative regime and prove 
exponential convergence towards the equilibrium.\end{abstract}



\thanks{C.S. has been supported by the APART research grant funded by
the Austrian Academy of Sciences.}
\maketitle

\begin{center}
Version: May 30th, 2007
\end{center}

\section{Introduction and main results}\label{secint}

In recent years the rigorous mathematical study of kinetic equations has been enlarged to a class of models which take into account also 
quantum effects, \cf \cite{Vi} for a general overview. 
This so-called \emph{quantum kinetic theory} can be seen as an attempt to incorporate certain properties of 
an underlying quantum systems into the framework of classical statistical mechanics. One might hope that these 
``hybrid models'' on the one hand allow for a somewhat simpler description of the particle dynamics, while maintaining, on the other hand, 
some purely quantum mechanical features such as generalized statistics for \emph{Bosons} and \emph{Fermions}. 
Clearly such models can only be justified 
in a \emph{semi-classical regime}, respectively, in situations where the transport properties of the 
particles are mainly governed by Newtonian mechanics. Indeed this point of view has already been adopted in the classical paper by 
Uehling and Uhlenbeck \cite{UU}, in which they derived their celebrated nonlinear Boltzmann type equation for quantum particles. 

Following their spirit most of the quantum kinetic models studied so far invoke nonlinear 
collision operators of Boltzmann type, 
see, \eg, \cite{Do, EMV, Lu1, Lu2}. 
For these kind of models, a focus of interest is on the long time behavior of their solutions, in particular 
the convergence towards steady states, which generalize the classical \emph{Maxwellian distribution}, 
\cf \cite{EMVe, Lu1, LuWe, NeSc}. 
Very often though, the simplified case of a spatially homogeneous gas is considered. 

In the present work, we shall also be interested in such kind of relaxation-to-equilibrium phenomena, in the 
\emph{spatially inhomogeneous} case. We shall not deal with a Boltzmann type equation, but rather 
study a \emph{nonlinear Fokker-Planck} type model (FP). 
More precisely, we consider 
\be \label{m1}
\partial_t f + p \cdot\nabla_x (f+\kappa f^2)=  \, \diverg_p  \left(\nabla_p f+p f 
(1+\kappa f)\right) ,
\ee
where, for any $t\geq 0$, $f=f(t,x,p)\geq 0$ denotes the \emph{particle distribution} on \emph{phase space} $ \Omega_x \times \R_p^d$. 
In what follows, the spatial domain is chosen to be $\Omega_x = \mathbb T^d$, the $d$-dimensional torus. 
This setting can be seen as a convenient and mathematically simpler 
replacement for the incorporation of confining potentials $V(x)$, needed to guarantee the existence of 
nontrivial steady states in the whole space. 
In \eqref{m1} we set $\kappa = -1$ for Fermions and $\kappa = 1$ for Bosons.  
For $\kappa =0$ equation \eqref{m1} simplifies to the classical \emph{linear}  
Fokker-Planck equation (or \emph{Kramer's equation}) on phase space.
For this linear model, the convergence to equilibrium has recently been studied in 
\cite{DeVi, HeNi, MoNe}, using several different approaches. 

The FP type model \eqref{m1} has been introduced in \cite{KaQu}, for classical 
particles obeying an exclusion principle. A formal derivation 
from a generalized Boltzmann equation for fermions and bosons is given in \cite{K95}, 
and from the Uehling-Uhlenbeck equation in \cite{RoKa}. 
Different physical applications can be found in \cite{Fr, K01, LKQ} dealing with, both, the spatially homogeneous 
as well as the inhomogeneous case (see also \cite{Fr1} and the references therein). 
More recently a similar but somewhat simpler FP type model has been proposed in \cite{SSC1, SSC2} 
to describe self-gravitating particles and the formation of Bose-Einstein condensates in a kinetic framework. 
The authors consider
\be \label{m2}
\partial_t f + p\cdot\nabla_x f=\, \diverg_p  \left(\nabla_p f+p f 
(1+\kappa f)\right),
\ee
where, in contrast to \eqref{m1}, only the diffusive part of the equation includes a nonlinearity. 
As we will see both equations however share the same steady states. 
In order to deal with both models at the same time we study from now on the following initial value problem
\be \label{model}
\left \{
\begin{aligned}
& \, \partial_t f + p \cdot\nabla_x (f+ \sigma \kappa f^2)=  \, \diverg_p  \left(\nabla_p f+p f 
(1+\kappa f)\right) ,\\
& \, f\big |_{t=0}=  f_0(x,p),
\end{aligned}
\right.
\ee
with $\sigma = 1$, or $\sigma =0$, corresponding to the case \eqref{m1} and \eqref{m2}, respectively. We note that 
the long time behavior of these models in the spatially homogeneous case (and for $d=1$) has been 
rigorously investigated quite recently in \cite{CRS} via an entropy-dissipation approach.
 
In what follows, the initial phase space distribution $f_0 \in \L^1(\mathbb T^d_x \times \R^d_p)$ is assumed to be normalized according to
\be\label{mass}
\iint_{\T^d\times \R^d} f_0(x,p) \, \D x \, \D p = M, 
\ee
for some given \emph{mass} $M >0$. This normalization is conserved by the evolution. 
Moreover in the fermionic case, \ie $\kappa =-1$, 
we require $f_0(x,p) < 1$, $\forall \, (x,p) 
\in \mathbb T^d \times \R^d $, as usual in the physics literature \cite{Fr}. 
In particular the latter is needed to define the associated quantum mechanical \emph{entropy functional}, \ie
$$
H[f] := \iint_{\T^d\times \R^d} \left( \frac{|p|^2}{2} f +    f \ln f 
- \kappa (1+\kappa f) \ln (1+ \kappa f) \right) \D x \, \D p,
$$
which obviously requires $f(t,x,p) < 1$, if $\kappa = -1$. 
It is now straightforward to verify that (independent of the particular choice of $\sigma$) 
the unique \emph{steady state} of \eqref{model} is given by
\begin{equation} \label{finfty} 
f_\infty  =\frac{1}{\exp{\left(\frac{|p|^2}{2}+ \theta \right)}-\kappa},
\end{equation}
where the constant $\theta $ is used to ensure that $f_\infty$ satisfies the 
mass constraint \eqref{mass}. 
In the bosonic case we require $\theta \in \R_+$, whereas in the fermionic case 
we can allow for $\theta \in \R$, \cf  \cite{CRS, EMV1} for more details. In the latter situation the distribution \eqref{finfty} 
is the well known \emph{Fermi-Dirac equilibrium distribution}. On the other hand, for 
$\kappa = 1$,  $f_\infty$ is the so-called \emph{regular Bose-Einstein distribution}. 
Finally, if $\kappa =0$, formula \eqref{finfty} simplifies to the classical Maxwellian, \ie 
$$
f^{\rm lin}_\infty  = \frac{M}{(2\pi)^{d/2}} \, e^{- |p|^2/2} , 
$$
where $\log (M(2\pi)^{-d/2}) = -\theta$. Note that in any case the equilibrium state is \emph{independent} of $x$ since 
we have chosen our spatial domain to be $\mathbb T^d$. 

In the bosonic case there is an additional difficulty, at least for $d\geq 3$, since 
one can \emph{not} associate an arbitrary large $M>0$ to the steady state. 
More precisely, the maximum amount of mass comprised by 
$f_\infty$ is determined via
$$
\iint_{\T^d\times \R^d} \frac{1}{e^{|p|^2/2}-1} \, \D x \, \D p =: M_{\rm crit} < \infty,
$$
\ie for $\theta = 0$. 
Due to mass conservation this induces a threshold on $M$. This problem, 
which does not appear in dimensions $d=1$ or $2$, has led to the introduction of 
more general bosonic steady states, where an additional $\delta$-distribution (appropriately normalized) 
is added to $f_\infty$, \cf \cite{EMV, EMVe}. 
This singular measure can then be interpreted as a so-called \emph{Bose-Einstein condensate} (BEC). 
The formation of a $\delta$-measure in finite or infinite time is a task of extensive research in 
quantum kinetic theory, see, \eg, \cite{BMP, EMVe}. 
For our nonlinear model though, including such generalized solutions on a rigorous mathematical level 
seems to be out of reach so far and we 
thus have to impose $\theta >0$ in the bosonic case. (Indeed, as we shall see below, we also require 
$\theta >0$ in the fermionic case, although for different and rather technical reasons.)

Our main task here is the description of the convergence for solutions of \eqref{model} towards the 
steady state \eqref{finfty}. To this end we shall conceptually follow the approach given in \cite{MoNe} where the trend 
to equilibrium is studied for a wide class of kinetic models \emph{close to equilibrium}. The difference in our case being mainly 
that we are dealing with \emph{local} nonlinearities which moreover 
are also allowed to enter in the transport part of the considered equation.
Mathematically speaking, the approach is based on the so-called \emph{hypo-coercivity} 
property of the linearized equation \cite{Vi1, Vi2}. 

We consequently \emph{linearize} the solution $f$ of \eqref{model} 
around the steady state $f_\infty$ in the form
\be\label{lin}
f = f_\infty + g \sqrt{\mu_\infty} .
\ee
where the new unknown $g(t,x,p)\in \R$ can be interpreted as a \emph{perturbation} of the equilibrium state such that 
$$
\iint_{\T^d\times \R^d}  g \sqrt{\mu_\infty} \, \D x \, \D p =0.
$$
In \eqref{lin} we use the additional (time-independent) \emph{scaling factor}
$$
\mu_\infty:= f_\infty +\kappa f^2_\infty ,
$$
which allows for an easier description in the functional framework given below. 
Plugging \eqref{lin} into \eqref{model} straightforward calculations formally yield the following equation for $g$
\begin{equation}\label{nlin}
\partial_t g + (1+ 2\sigma \kappa f_\infty ) p\cdot\nabla_x g  =  L(g)+ Q(g)\ .
\end{equation}
Above the \emph{linearized collision operator} $L$ is given by
\be\label{lincollop}
\begin{aligned}
L (g)= & \ \frac{1}{\sqrt{\mu_\infty}} \,\diverg_p \big(\nabla_p \,( g \sqrt{\mu_\infty} )+ p\, \eta_\infty g \, \sqrt{\mu_\infty} \, \big)\\
= & \ \Delta_p g+g\left(\frac{d}{2} \, \eta_\infty - |p|^2\left(\frac{1}{4}+2 \kappa \mu_\infty \right)\right) ,
\end{aligned}
\ee
where we use the short hand notation $\eta_\infty :=1+2 \kappa f_\infty $. 
The \emph{quadratic remainder} $ Q$ is 
\be \label{quadratic}
Q ( g)= \frac{\kappa}{\sqrt{\mu_\infty}} \, \big (\diverg_p \left(p \, \mu_\infty g^2\right) -  
\sigma \mu_\infty p\cdot \nabla_x (g^2)  \big).
\end{equation}

The main result of our work is as follows.
\begin{theorem}\label{th1} 
Let $f_0$ be of the form 
$$
0\leq f_0 = f_\infty + g_0 \sqrt{\mu_\infty}\,,
$$
with $\theta >0$. Moreover if $\kappa = \sigma = 1$, \ie the bosonic case with nonlinear transport, 
assume that in addition $\theta >  \theta^*$, for a certain $\theta^* >0$. 

For $k\in \mathbb N $ with $k> 1+d/2$ there exists an $\epsilon_0 >0$, such that for all $f_0$ with 
$\norm{\, g_0 \, }_{H^k}\leq\epsilon_0$, the equation \eqref{model} 
admits a unique solution 
$0\leq f\in C([0,\infty); H^k(\mathbb T^d_x\times \R^d_p))$. Moreover
\begin{equation*}
\norm{\, \mu_\infty^{-1/2} \big(f(t)-f_\infty \big)}_{H^k}\leq  C(\epsilon_0)\, e^{-\tau t},
\end{equation*}
where $C(\epsilon_0)$ and $\tau$ are positive constants.
\end{theorem} 

First we collect several preliminary 
results in the Section \ref{seclin}. The proof of Theorem \ref{th1} is then given in Section \ref{secproof}.

\begin{remark} \label{remark}\mbox{}
\begin{itemize}
\item For Bosons it 
is crucial to avoid the possible formation of a BEC by imposing $\theta>0$ since formal calculations given in \cite{SSC2} indicate that 
an analogous theorem can not hold if $M > M_{\rm crit}$, see also \cite{EMVe}. In the fermionic case the 
reason to impose $\theta >0$ is to guarantee $\eta_\infty>0$, which we will make use of several times. 
It might be possible to overcome this restriction for fermions by using a different approach, see \cite{CRS}. 
\item The additional requirement $\theta > \theta^*$ is \emph{not} needed for the model \eqref{m2} where only a nonlinear diffusion 
operator is present. The reason for this constraint when dealing with \eqref{m1} 
is that we have to maintain a fundamental regularizing property in $x$ of the transport part, \cf 
Section \ref{secproof} for more details. 
Since $\theta^*$ is then determined by a 
transcendental equation, we do not give an exact value for $\theta^*$ but only perform numerical experiments which 
indicate that $\theta^* \approx 0.451$.
\item As already discussed in \cite{MoNe} one could also 
take into account \emph{self-consistent} potentials, which stem from 
a coupling to Poisson's equation. The latter case might be particularly interesting in semiconductor modeling, where 
the fermionic FP type equation could be used to describe the dynamical behavior of charge carries obeying the 
``physically correct'' equilibrium statistics. 
\end{itemize}
\end{remark}
\begin{corollary} Under the same assumptions as above 
$$
H[f(t)] - H[f_\infty] \leq C e^{-\tau t},
$$
\ie we have exponential decay in relative entropy.
\end{corollary}

\begin{proof} Inserting $f=f_\infty + \sqrt{\mu_\infty} g$ into $H[f]$, having in mind $\iint \sqrt{\mu_\infty} g \, \D x \, \D p=0$, we
perform a Taylor expansion around the steady state $f_\infty$, and finally use Theorem \ref{th1}. This yields the assertion of the 
corollary .
\end{proof}

\section{Study of the linearized collision operator}\label{seclin}

We shall now derive several properties of the linearized collision operator $L$ to be used in the proof of the main result. 
First note that $L$ is self adjoint on $L^2(\R^d_p)$ and, 
by partial integration, one obtains
\begin{equation}\label{collsym}
{\langle L(g),g \rangle}_{L^2(\R_p^d)}=-\int_{\R^d}\left|\, \nabla_p g +\frac{p}{2}\, \eta_\infty g \, \right|^2\ud p=-\int_{\R^d}\left|\nabla_p\left(\frac{g}{\sqrt{\mu_\infty}} \right)\right|^2\mu_\infty \ud p\ .
\end{equation}
Thus the kernel of the non-positive operator $L$ is given by
$$
\mbox{Ker}(L)=\mbox{span}\{ \sqrt{\mu_\infty} \, \}.
$$
Let us define the orthogonal projection in $L^2(\R^d_p)$ onto this kernel via 
$$
\Pi (f):= \left( \frac{ 1}{\rho_\infty} \int_{\R^d} f \, \sqrt{\mu_\infty} \, \ud p\right) \sqrt{\mu_\infty}\,, 
$$
where we set 
$$
\rho_\infty = \left(\int_{\R^d} \mu_\infty  \, \ud p \right) >0
$$
for reasons of normalization. Note that this is only a projection in the momentum variable $p \in \R^d$. 
Motivated by \eqref{collsym} we introduce the following weighted space
$$
\Lambda_p:=\big \{ f \in L^2(\R_p^d):\ \norm{ f }_{\Lambda_p} < \infty \big \}, 
$$
where
$$
\norm{ f }_{\Lambda_p}^2 : = \norm{\nabla_p f}_{L_p^2}^2+ \norm{ p \, \eta_\infty \, f  }_{L_p^2}^2 ,
$$
Here, and in what follows, we write $L^2_p\equiv L^2(\R^d_p)$ for simplicity. Moreover we denote by 
$$
\norm { f  }_\Lambda := \norm{ \,\norm{  f }_{\Lambda_p} }_{L^2(\mathbb T^d_x)},
$$
the induced norm on phase space. 
Obviously the $\Lambda_p$-norm controls the $L^2_p$-norm for $\kappa$ nonnegative. In the fermionic case ($\kappa = -1$) however  
this is not true in general since $\eta_\infty$ may change sign. For our functional approach the control of $L^2_p$ via 
$\Lambda_p$ is crucial and thus we have to guarantee that $\eta_\infty >0$ by assumption. This implies that for $\kappa = -1$ 
we need to impose $f_\infty < 1/2$, $\forall  \, p \in \R^d$, or equivalently $\theta >0$. This 
certainly is more restrictive than the usual bound \ie $f_\infty < 1$, used in the physics literature. 
In summary we require $\theta>0$ in the bosonic case to prevent BEC and in the 
fermionic case to ensure that $\eta_\infty$ is globally bounded away from zero.

First we obtain a Poincar\'{e} inequality for the steady state of the linearized model.
\begin{lemma} If $\theta >0$ then the (strictly positive and normalized) measure $\mu_\infty/\rho_\infty$ satisfies a Poincar\'{e} 
inequality on $\R^d_p$, \ie
$$
\int_{\R^d} g^2 \mu_\infty \, \D p - \frac{1}{\rho_\infty}\left(\int_{\R^d} g \mu_\infty \, \D p \right)^2 \leq C_{\rm p} 
\int_{\R^d} |\nabla_p g|^2 \mu_\infty \, \D p, \quad C_{\rm p}>0. 
$$
\end{lemma}
Note that in the fermionic case, $\kappa = -1$, this lemma indeed holds more generally for any $\theta \in \R$. 
\begin{proof}  
Let $\mathcal A $ be defined by $\mathcal A = -\ln \mu_\infty$. Inequality $(3.6)$ in \cite{AMTU} shows that if $\mathcal{A}$ is uniformly convex then a Poincar\'{e} inequality with measure $\mu_\infty/\rho_\infty$ (where we divide by $\rho_\infty$ for normalization purposes) holds. Moreover Theorem~3.3. in \cite{AMTU} allows to take into account also $L^\infty$ perturbations of $\mathcal{A}$.
Thus we wish to decompose
$$
\mathcal A  =  \mathcal A_1  +  \mathcal A_2, 
$$ 
where $\mathcal A_1 (p)$ is uniformly convex and $\mathcal A_2(p)$ is a $L^\infty$ perturbation. To this end note that $\mathcal{A}$ is given by
\begin{align*}
\mathcal{A}= & \ -\log\left(\frac{e^{|p|^2/2+\theta}}{(e^{|p|^2/2+\theta}-\kappa)^2}\right) =  \frac{|p|^2}{2}+\theta-2\log \left(\frac{e^{|p|^2/2+\theta}}{e^{|p|^2/2+\theta}-\kappa}\right)\\
= & \ \frac{|p|^2}{2}+\theta+2\log\left(1-\frac{1}{\kappa e^{|p|^2/2+\theta}}\right)\ .
\end{align*}
We now pick $\mathcal{A}_1=|p|^2/2$ and $\mathcal{A}_2$ to be the rest of the terms appearing on the right hand side. Then $\mathcal{A}_2$ is bounded if $\kappa = -1$, or if $\kappa =1$ and $\theta >0$ which concludes the proof.
\end{proof}

With the above lemma in hand we can now establish the \emph{coercivity} of the linearized collision operator. From now on we will denote by $K, K_1,\dots$ generic (positive) constants, to be used several times in different proofs, 
whereas the constants $C_1, C_2, \dots$ will be fixed throughout the work in order to be able to track their appearance.

\begin{lemma}\label{coerc}
For $\theta>0$ there exists a $\lambda >0$ such that 
$$
\langle L(g),g\rangle_{L_p^2}\leq -\lambda \norm{g-\Pi(g)}_{\Lambda_p}^2, \quad \forall \, g\in \Lambda_p\,.
$$
\end{lemma} 
The coercivity property (in $p\in \R^d$) of the operator $L$ is indeed an essential requirement 
to establish our main result. 
\begin{proof}
We start with
\begin{align*}
{\langle L(g),g\rangle}_{L_p^2}= & \ -  \rho_\infty\int\left|\nabla_p\left(\frac{g}{\sqrt{\mu_\infty}}\right)\right|^2\frac{\mu_\infty}{\rho_\infty}\, \ud p\\
 \leq & \ - C_{\rm p} \left (\int g^2\ud p-\frac{1}{\rho_\infty}\left(\int\sqrt{\mu_\infty}\, g\ud p\right)^2\right) \\
 \leq & \ - K_1 \int\left(g-\frac{\sqrt{\mu_\infty}}{\rho_\infty}\int\sqrt{\mu_\infty} \, g\ud p\right)^2\ud p=
 \, - K_1 \norm{ g-\Pi(g)}_{L^2_p}^2 ,
\end{align*}
for some $K_1>0$, 
where we used the fact that the measure $\mu_\infty/\rho_\infty$ satisfies a Poincar\'{e} inequality due to the previous lemma.
Now  to improve on the $L_p^2$-norm we use
\begin{align*}
{\langle L(g),g\rangle}_{L_p^2}= & - \ \int\left(\nabla_p(g-\Pi(g))+ \frac{p}{2}\eta_\infty(g-\Pi(g))\right)^2\ud p \\
\leq & \  - K_2 \norm{g-\Pi(g)}_{\Lambda_p}^2 + K_3 \norm{g-\Pi(g)}_{L_p^2}^2\ .
\end{align*}
Adding the two inequalities above (multiplied by appropriate constants) finishes the proof.
\end{proof}

We get a similar result for the derivatives w.r.t. $p \in \R^d$.

\begin{lemma}\label{coercdp}
Let $\theta >0$, then there exist positive constants $C_1$ and $C_2$, such that for any $g \in L^2_p$ with 
$\nabla_p g\in\Lambda_p$
$$
{\langle\nabla_p L(g),\nabla_p g\rangle}_{L_p^2}\leq - C_1 \norm{\nabla_p g}_{\Lambda_p}^2+ C_2 \norm{g }_{L_p^2}^2.
$$
\end{lemma}
\begin{proof}
A lengthy calculation yields 
\begin{align*}
{\langle\nabla_p L(g),\nabla_p g\rangle}_{L_p^2}=  & \, \int - (\Delta_p g)^2 + |\nabla_p g|^2 \left(\frac{d}{2} \eta_\infty - 
\frac{|p|^2}{4} - 2 \kappa \mu_\infty |p|^2\right)\, \D p\\
& \,  + \int \frac{g^2}{2} \left( d^2\kappa\mu_\infty+\frac{d}{2}+4\kappa d\mu_\infty\right)\,\ud p\,\\
& \, - \frac{\kappa}{2} \int g^2 |p|^2\mu_\infty\left((d+10)\eta_\infty-2|p|^2\mu_\infty\eta_\infty\right)\,\ud p.
\end{align*}
The last integral on the r.h.s. is dominated by the $L_p^2$-norm, since 
$\mu_\infty$ decays exponentially fast as $|p|\to \infty$. 
We also have that $1/4 + 2 \kappa \mu_\infty \geq C>0 $, $\forall \, p \in \R^d$. This obviously holds true for the bosonic 
case but is also guaranteed in the fermionic situation where $f_\infty < 1/2$. Thus we can estimate 
\begin{align*}
{\langle\nabla_p L(g),\nabla_p g\rangle}_{L_p^2} \leq - C_1 \norm{\nabla_p g}_{\Lambda_p}^2 + 
\frac{d}{2} \int \eta_\infty |\nabla_p g|^2 \, \D p + C_2 \norm{g }_{L_p^2}^2
\end{align*}
and a classical interpolation argument applied to the second term on the r.h.s. 
yields the assertion of the lemma (with different constants $C_1$, $C_2$).
\end{proof}

Finally we need the following technical lemma.

\begin{lemma} \label{split}
For $g,h \in \Lambda_p$ it holds that
$$
\langle L(h),g\rangle_{L_p^2}\leq C_3 \norm{ g}_{\Lambda_p} \norm{h}_{\Lambda_p}, \quad C_3>0.
$$
\end{lemma}
The combination of Lemma \ref{coerc} and Lemma \ref{split} induces the particular choice of $\Lambda_p$ and its corresponding norm. 
\begin{proof}
We first note that
\begin{align*}
& \norm{ g}_{\Lambda_p} \norm{h}_{\Lambda_p} = \left(\norm{\nabla_p g}_{L_p^2}^2+ \norm{ p \, \eta_\infty \, g  }_{L_p^2}^2 \right)^{\frac{1}{2}}
\left(\norm{\nabla_p h}_{L_p^2}^2+ \norm{ p \, \eta_\infty \, h  }_{L_p^2}^2 \right)^{\frac{1}{2}}\\
& \ \geq \frac{1}{C_3} \, \left(  \norm{\nabla_p g}_{L_p^2}^2 + \int (1+|p|^2) \, \eta^2_\infty g^2 \, \D p \right)^{\frac{1}{2}} 
\left(  \norm{\nabla_p h}_{L_p^2}^2 + \int (1+|p|^2) \, \eta^2_\infty h^2 \, \D p \right)^{\frac{1}{2}}, 
\end{align*}
for some $C_3> 0$ since the $\Lambda_p$-norm dominates the $L^2_p$-norm. Using the following simple algebraic estimate 
\begin{equation*}
\big((a^2+b^2)(c^2+d^2)\big)^{\frac{1}{2}}\geq ac+bd\ ,
\end{equation*}
(with $a,b,c,d \in \R$) we further obtain 
\begin{equation*}
C_3 \norm{ g}_{\Lambda_p} \norm{h}_{\Lambda_p} \geq \norm{\nabla_p g}_{L_p^2}\norm{\nabla_p h}_{L_p^2}+\left( \int (1+|p|^2) \, \eta^2_\infty g^2 \, \D p \int (1+|p|^2) \, \eta^2_\infty h^2 \, \D p \right)^{\frac{1}{2}}.
\end{equation*}
The proof then follows by applying the Cauchy-Schwarz inequality to both terms on the right hand side 
and integrating by parts in the first one.
\end{proof}
 
\section{Convergence for the linear model and proof of Theorem \ref{th1}}\label{secproof} 

Now we are able to establish the long time asymptotics for the linearized equation, which eventually 
will be translated also to the nonlinear model \eqref{model} in the perturbative setting.
\begin{proposition}
Consider the linearized Fokker Planck type equation
\begin{equation}\label{linearized}
\partial_t g + (1+ 2\sigma \kappa f_\infty ) p\cdot\nabla_x g  =  L(g)\ ,
\end{equation}
with $L$ given by \eqref{lincollop} and $\theta >0$. Moreover if $\kappa = \sigma = 1$, assume 
in addition that $\theta > \theta^*$. 
Let the initial data $g_0\in H^k(\mathbb{T}^d\times \mathbb{R}^d)$, for $k \in \N$. Then
the solution $g(t)$ exists globally in time and 
\begin{equation*}
\norm{\, g(t) -g_\infty}_{H^k}\leq C\, e^{-\tau t}\ , \quad \mbox{ with $C=C({\|g_0\|}_{H^k})$, $\tau >0$,}
\end{equation*}
where the global equilibrium $g_\infty$ is given by
\begin{equation*}
g_\infty=\left( \frac{ 1}{\rho_\infty} \iint_{\T^d\times \R^d} g_0 \, \sqrt{\mu_\infty} \, \ud p\, \ud x\right) \sqrt{\mu_\infty}\ .
\end{equation*}
\end{proposition}
We want to remind the reader of the discussion of the role of $\theta^*$ in Remark~\ref{remark}.
\begin{proof}
For the proof we proceed similarly to \cite{MoNe}. We will sketch the main ideas and 
stress the differences which occur due to the changes in the transport operator.
Note that, since the equation is linear, we can w.r.o.g. consider the case where 
$g_\infty \equiv 0$. This can always 
be achieved by subtracting initially the projection onto the global equilibrium, 
\ie by considering initial data $\widetilde g_0=g_0-g_\infty$.

We start with $k=1$. 
The main idea of the proof is to study the time evolution of a combination of derivatives w.r.t. $x$ and $p$. 
More precisely we consider the following functional
\begin{equation*}
\mathcal F [g(t)]:=\alpha\norm{g}^2+\beta \norm{\nabla_x g}^2+\gamma\norm{\nabla_p g}^2+\delta\left<\nabla_xg, \nabla_pg\right>\ ,
\end{equation*}
where $\norm {\, \cdot \,}$ denotes the standard norm on $L^2(\mathbb T_x^d \times \R_p^d)$ and $
\alpha, \beta, \gamma, \delta$ are some positive constants.
We note that $\delta$ has to be small enough in 
comparison to $\beta$ and $\gamma$ such that $\mathcal F$ is positive and controlled from above and 
below by the square of the usual $H^1(\mathbb T_x^d \times \R_p^d)$-norm of $g$. On the 
other hand $\delta$ has to be strictly positive, since we need it in order to close the argument (see below). We 
aim to prove that 
\begin{equation}
\label{mainest}
\dt \ \mathcal{F}[g(t)]\leq-\tilde C\left(\norm{ g} _\Lambda^2+ \norm{\nabla_{x,p} \, g} _\Lambda^2\right), \quad \tilde C >0.
\end{equation} 

To this end we calculate the time derivatives of the various summands in $\mathcal{F}$. 
First, for the $L^2$-norm we have
\begin{equation}\label{normg}
\dt \, \norm{g}^2=2\left<L(g),g)\right>\leq - 2 \lambda  \norm{g-\Pi(g)}_{\Lambda}^2,
\end{equation}
where we have used that the transport part does not contribute due to its divergence form (which is a straightforward calculation) and 
the assertion of Lemma \ref{coerc}. 

Next the spatial derivatives evolve according to 
\begin{equation} \label{xder}
\dt \, \norm{\nabla_x g}^2=2\left<\nabla_xL(g),\nabla_xg\right>\leq - 2\lambda \norm{\nabla _x g-\Pi( \nabla _x g)}_{\Lambda}^2,
\end{equation}
where we again used that the contribution from the transport term vanishes and the fact that $L$ commutes with $\nabla_x$, 
thus allowing us to apply Lemma \ref{coerc} also on $\nabla_x g$. 

For the derivatives w.r.t. $p$ we get some additional terms by the coefficient of the transport operator
\begin{align*}
\dt \, \norm{ \nabla_p g}^2=  & \ 2\left<\nabla_pL(g),\nabla_pg\right> 
+\left<\left(4\kappa\sigma\mu_\infty|p|^2-2\eta_\infty(1+2\kappa\sigma f_\infty)\right)\nabla_xg,\nabla_pg\right>\ \\
\leq & \  2\left<\nabla_pL(g),\nabla_pg\right> + K \iint |\nabla_p g| |\nabla_x g| \, \D p \, \D x.
\end{align*}
Here we have used that $f_\infty$ as well as $|p|^2\mu_\infty$ 
are uniformly bounded in $L^\infty$. 
Now the first term on the r.h.s. is estimated by Lemma \ref{coercdp}, which yields some damping 
for $\nabla_p g$ in the $\Lambda$-norm, \ie
\begin{equation} \label{zwischen}
\begin{aligned}
\dt \, \norm{ \nabla_p g}^2\leq & \, 
- 2C_1 \norm{\nabla_p g}_{\Lambda}^2+ 2C_2 \norm{g }^2 + K \iint |\nabla_p g| |\nabla_x g| \, \D p \, \D x \\
\leq & \, -  K_1 \norm{\nabla_p g}_{\Lambda}^2+ K_2 \norm{\nabla_x g}^2 + 2C_2 \norm{g }^2,
\end{aligned}
\end{equation}
with $K_1, K_{2} >0$. In order to deal with the term $\propto \norm{g}^2$, we split 
$$
g=(g-\Pi(g))+\Pi(g)
$$ 
and estimate 
\begin{align*}
\norm{g}^2 \leq  \norm{ g-\Pi(g)}^2+ \norm{ \Pi(g)}^2 \leq
 \norm{ g-\Pi(g)}^2+ C_{\mathbb T} \norm{ \nabla_x g}^2, \quad C_{\mathbb T}>0.
\end{align*}
In the second step, we used the classical Poincar\'{e} inequality w.r.t 
$x\in \mathbb T^d$ and the fact that $\Pi(g)$ has zero mean on the torus, since $\iint g_\infty\ud p\ud x = 0$. 
Together with \eqref{zwischen}, this yields
\begin{equation}\label{pder}
\dt \, \norm{ \nabla_p g}^2\leq - K_1 \norm{\nabla_p g}_{\Lambda}^2+2C_2\norm{ g-\Pi(g)}^2+K_3 \norm {\nabla_x g} ^2 ,
\end{equation} 
where $K_3= K_2+2C_2C_{\mathbb{T}}$. Note that this is an 
improvement to \eqref{zwischen}, since 
the term $\|\, g-\Pi(g) \, \|^2$ can be controlled by adjusting $\alpha$ and having in mind \eqref{normg}.
 
Finally we look at the mixed derivatives w.r.t. $x$ and $p$, which evolve according to
\begin{equation}\label{mixedder}
\dt \, \left<\nabla_x g,\nabla_p g\right>=2\left<L(\nabla_xg),\nabla_pg\right>-\left<\nabla_xg,
\big(1+2\sigma\kappa f_\infty-2\sigma\kappa\mu_\infty|p|^2)\big )\nabla_xg\right>.
\end{equation}
For the first term on the right hand side we invoke Lemma \ref{split}, which together with the 
Cauchy-Schwarz inequality in $x$ implies 
$$
\left<L(\nabla_xg),\nabla_pg\right> \leq C_3 \eta \norm{\nabla _x g-\Pi( \nabla _x g)}_{\Lambda}^2 +  C_3 \eta^{-1} \norm{\nabla_p g}^2_\Lambda, \quad 
\forall \, \eta >0.
$$
The second term on the right hand side of \eqref{mixedder}, which stems from the 
transport part, generates a damping for $\norm{\nabla_x g}$ 
(which the operator $L$ can not provide since it only acts 
in $p$), provided that 
\be \label{addcon}
1+2\kappa \sigma f_\infty-2\kappa \sigma \mu_\infty|p|^2\geq K_4 >0 .
\ee
Assuming for the moment that \eqref{addcon} is true we obtain
\begin{equation*}
\dt \, \left<\nabla_x g,\nabla_p g\right>\leq 2C_3\eta \norm{\nabla _x g-\Pi( \nabla _x g)}_{\Lambda}^2+2 C_3\eta^{-1} \norm{\nabla_p g}^2_\Lambda-K_4\norm{\nabla_x g }^2\ .
\end{equation*}
In summary this estimate, together with \eqref{normg}, \eqref{xder}, and \eqref{pder}, yields
\begin{align*}
\dt \ \mathcal{F}[g(t)]\leq & \, - 2 (\lambda \alpha - C_2 \gamma)\norm{ g-\Pi( g)}_{\Lambda}^2 
- 2 (\beta \lambda -C_3 \eta \delta) \norm{\nabla _x g-\Pi( \nabla _x g)}_{\Lambda}^2 \\
& \, -  (K_1 \gamma -  2 C_3 \delta \eta^{-1})\norm{\nabla_p g }_\Lambda^2 -(K_4 \delta - K_3 \gamma) \norm{\nabla_x g }^2\ .
\end{align*}
It remains to find coefficients $\alpha,\beta,\gamma, \delta, \eta$ in $\mathcal{F}$, such such that all ``bad'' terms in the above 
given estimates (\ie those which come with the wrong or without sign) can be controlled 
and the differential inequality \eqref{mainest} holds true. 
This can be done analogously to Step 4 in the Proof of Theorem 1.1 in \cite{MoNe} and we therefore will not elaborate further on it. 
Note that the functional $\mathcal F$ then clearly 
induces a new norm on phase space, equivalent to $H^1(\R^d\times \mathbb T^d)$, via
$\norm{\, g \, }_{\mathcal H^1}^2 := \mathcal F[\, g \, ].$

To retain the (fundamental) damping property in the spatial derivatives coming from the evolution of the mixed term it 
remains to show that the constraint \eqref{addcon} holds true. 
We denote
\begin{equation*}
\Psi_\kappa (p\, ; \theta):= 1+2\kappa f_\infty-2\kappa\mu_\infty|p|^2 \ .
\end{equation*}
If $\kappa = -1$ and since $\theta>0$, 
it is clear that $\Psi_{-1}(p\, ; \theta)\geq K_4>0$ because $\eta_\infty \geq 0$ in this case. 
For Bosons, \ie $\kappa = 1$, however the situation is more difficult. 
Note that 
$$
\lim _{|p|\to \infty} \Psi_1(p\, ; \theta) = \lim _{|p|\to \infty} (\eta_\infty-2\mu_\infty|p|^2) = 1, \quad \forall \, \theta >0,
$$ 
and thus by continuity it is enough to make sure that $\Psi_1(p\, ; \theta) \not = 0$, $\forall \, p \in \R^d, \theta>\theta^*$.
Straightforward calculations show that $\Psi_1(p\, ; \theta)$ can only be zero if
\begin{equation*}
\left(e^{|p|^2/2+\theta}\right)^2-2|p|^2e^{|p|^2/2+\theta}-1=0\ ,
\end{equation*}
which implies  
\begin{equation*}
e^{|p|^2/2+\theta}=|p|^2+\sqrt{|p|^4+1}\ .
\end{equation*}
Obviously this equality can not be true for $\theta $ larger than some critical value $\theta^*$. 
Numerical experiments suggest that this critical value is approximately $\theta^*\approx 0.451$. 
In summary one obtains the final estimate  \eqref{mainest}, which finishes the proof for $k=1$.

To proceed to higher order estimates in the Sobolev index $k \in \mathbb N$, we observe 
that the proof of Lemma \eqref{coercdp} can be 
generalized in a straightforward way to obtain 
\begin{equation*}
{\langle \partial_{x_\ell} \partial_{p_j} L(g),\partial_{x_\ell} \partial_{p_j} g\rangle}_{L_p^2}\leq - 
C_{1,k} \norm{\, \partial_{x_\ell} \partial_{p_j} g}_{\Lambda_p}^2+ C_{2,k} \norm{g }_{H^{k-1}_p}^2,
\end{equation*}
for any multi-indices $j, \ell$ such that $k=|j|+|\ell|$, $|j|\geq 1$. An induction argument in $k \in \mathbb N$, 
similar to the one given in \cite[Theorem 3.1]{MoNe} then yields the corresponding statement in $H^k$. 
There are no additional problems due to 
the coefficient in the transport operator, since the terms containing the highest order derivatives of $g$ 
can be treated as above and the lower order terms, which contain derivatives of $\Psi_\kappa(p \,;\theta)$, 
can be handled by interpolation since $\partial_{p_j}\Psi_\kappa(p \,;\theta) \in L^\infty$. 
This concludes the proof. 
\end{proof} 

Now we apply the result for the 
linearized equation to the nonlinear problem.

\begin{proof}[Proof of Theorem \ref{th1}] 

We have to show that the quadratic nonlinearity does not change the estimates obtained for the 
linearized equation, as long as the deviation from the equilibrium is small. 
The function $g = (f-f_\infty) \mu_\infty^{-1/2}$ solves \eqref{nlin} from which we 
deduce 
\begin{align*}
\frac{\D}{\D t} \norm{ g }_{\mathcal H^k}^2 = & \ 2 \left< T g, g \right >_{\mathcal H^k} + 2 \left< Q( g), g \right >_{\mathcal H^k},
\end{align*}
where $T:= L -(1+ 2\sigma \kappa f_\infty ) p\cdot\nabla_x$ and $L, Q$ 
are given in \eqref{lincollop}, \eqref{quadratic}, respectively. From the proof of 
Proposition \ref{linearized} we know that 
\begin{align*}
\left< T g, g \right >_{\mathcal H^k} \leq \, - \,\tilde  C \left(\sum_{|j|+|\ell |\leq k} \norm{ \, \partial_{x_\ell} \partial_{p_j} g}^2_\Lambda\right),
\end{align*}
where $\tilde C$ is, as in \eqref{mainest}. Thus, if we can prove the following property for the nonlinear part
\begin{align}\label{qest}
\left< Q( g), g \right >_{ H^k}\leq C_Q \norm{ g }_{ H^k}^2 \left(\sum_{|j|+|\ell |\leq k} \norm{ \, \partial_{x_\ell} \partial_{p_j} g}_\Lambda\right),
\end{align}
it follows, since $\norm{  \, \cdot \, }_{\mathcal H^k} \simeq \norm{ \, \cdot\, }_{H^k}$, that 
\begin{multline*}
\frac{\D}{\D t} \norm{ \, g \,}_{\mathcal H^k}^2 \leq  
- \,2\,\tilde C \left(\sum_{|j|+|\ell |\leq k} \norm{ \, \partial_{x_\ell} \partial_{p_j} g}^2_\Lambda\right) +
C_Q \norm{ g }_{ H^k}^2 \left(\sum_{|j|+|\ell |\leq k} \norm{ \, \partial_{x_\ell} \partial_{p_j} g}_\Lambda\right) \\
\leq  - 2\,\tilde C \left(\sum_{|j|+|\ell |\leq k} \norm{ \, \partial_{x_\ell} \partial_{p_j} g}^2_\Lambda\right) +
\epsilon C_{Q}\norm{ g }_{ H^k}^2 + C_\epsilon C_Q \norm{ g }_{ H^k}^2 \left(\sum_{|j|+|\ell |\leq k} \norm{ \, \partial_{x_\ell} \partial_{p_j} g}_\Lambda^2\right).
\end{multline*}
Now choosing $\epsilon$ small enough, such that $\epsilon C_{Q} \leq \tilde C$, we derive
\begin{equation*}
\frac{\D}{\D t} \norm{ \, g \,}_{\mathcal H^k}^2 \leq - \tilde C \left(\sum_{|j|+|\ell |\leq k} \norm{ \, \partial_{x_\ell} \partial_{p_j} g}^2_\Lambda\right) + C^{*} \norm{ g }_{ H^k}^2 \left(\sum_{|j|+|\ell |\leq k} \norm{ \, \partial_{x_\ell} \partial_{p_j} g}_\Lambda^2\right).
\end{equation*}
This concludes the proof of Theorem \ref{th1} by maximum principle as long as $\|g_0\|_{H^k}$ is sufficiently small.

In order to prove \eqref{qest}, we recall that $Q(g)$ is given by 
\begin{align*}
Q ( g)=  \frac{\kappa}{\sqrt{\mu_\infty}} \, \big(\diverg_p (p \, \mu_\infty g^2) \big) - \frac{\kappa}{\sqrt{\mu_\infty}}\,   
\big (\sigma \mu_\infty p\cdot \nabla_x (g^2)  \big) \equiv  Q_1(g) - Q_2(g),
\end{align*}
and note that $\mu_\infty^{-1/2} \partial_{x_\ell} \partial_{p_j} (p \, \mu_\infty) \in L^\infty(\mathbb T^d_x \times \R^d_p)$, 
for all multi-indeces $\ell, j \in \N^d$. 
We shall now treat $Q_1(g)$ and $Q_2(g)$ separately, using that $H^k(\R^d) \subset L^\infty(\R^d)$, 
for $k > d/2$, together with Leibniz' formula to differentiate $Q(g)$. It is relatively easy to see that 
\eqref{qest} holds for $Q_1(g)$, by estimates in the spirit of \eqref{estsob} below, since the $\Lambda$ norm incorporates an 
additional derivate w.r.t $p \in \R^d$. The estimate for $\left<Q_2(g), g \right>_{H^k}$ is more complicated 
though, since $Q_2$ contains a derivate w.r.t. $x\in \mathbb T^d$ which is not taken into account for 
by the $\Lambda$-norm. 

Thus we have to estimate
$$
\left< Q_2 (g), g \right >_{ H^k}=-\kappa\sigma\,\sum_{|j|+|\ell |\leq k}\left<\partial_{x_\ell} \partial_{p_j}
\left(\sqrt{\mu_\infty}\, p\cdot \nabla_x (g^2)\right), \,  \partial_{x_\ell} \partial_{p_j} g \right>_{L^2}.
$$
Since $\partial_{x_\ell} \partial_{p_j} \sqrt{ \mu_\infty} \in L^\infty(\mathbb T^d_x \times \R^d_p)$ the highest order terms are of the form
$$
\left< \sqrt{\mu_\infty}\, p\cdot \nabla_x \partial_{x_\ell} \partial_{p_j} (g^2), \,  \partial_{x_\ell} \partial_{p_j} g \right>_{L^2}, \quad 
|j|+|\ell |= k\, .
$$
Moreover, because of the additional derivative w.r.t. $p$ in the $\Lambda$-norm, the most problematic terms are those where 
$|\ell|=k$. Denoting $\partial_\ell \equiv \partial_{x_\ell}$ we compute
\begin{align*}
\left< \sqrt{\mu_\infty}\, p\cdot \nabla_x \partial_{\ell}  (g^2), \,  \partial_{\ell}  g \right>_{L^2} = 
& \ 2 \left< \sqrt{\mu_\infty}\, g p\cdot \nabla_x \partial_{\ell} g, \,  \partial_{\ell}  g \right>_{L^2} \\ 
& \, + \sum_{ i = 1}^d \big< \sqrt{\mu_\infty}\, p_i \sum_{{0\leq r \leq \ell + \delta_i}\atop{0< |r| < |\ell| +1}} 
\begin{pmatrix}
\ell + \delta_i \\
r
\end{pmatrix} 
\partial_{r}  g \, \partial_{\ell + \delta_i - r} \, g , \,  \partial_{\ell}  g \big>_{L^2},
\end{align*}
where $\delta_i$ denotes the $i\,$-th standard basis vector in $\R^d$. Using divergence theorem, we obtain 
\begin{align*}
\left< \sqrt{\mu_\infty}\, p\cdot \nabla_x \partial_{\ell}  (g^2), \,  \partial_{\ell}  g \right>_{L^2} =
& \ - \left< \sqrt{\mu_\infty}\, (p\cdot \nabla_x g ) \, \partial_{\ell} g, \,  \partial_{\ell}  g \right>_{L^2}\\
& \, + 
\sum_{ i = 1}^d \big< \sqrt{\mu_\infty}\, p_i \sum_{{0\leq r \leq \ell + \delta_i}\atop{1\leq |r| < |\ell| }} 
\begin{pmatrix}
\ell + \delta_i \\
r
\end{pmatrix} 
\partial_{r}  g \, \partial_{\ell + \delta_i - r} \, g , \,  \partial_{\ell}  g \big>_{L^2} .
\end{align*}
To estimate the first term on the right hand side we note that (remember $|\ell|=k$)
\begin{align}\label{estsob} 
\left<\sqrt{\mu_\infty}\,  (p\cdot \nabla_x g ) \, \partial_{\ell} g, \,  \partial_{\ell}  g \right>_{L^2} \leq K_1
\norm{\, \partial_\ell g \, }_{L^2}^2 \, \norm{\, \nabla_x g \,  }_{L^\infty} 
\leq K_2 \norm{\,  g \, }_{H^k}^2 \, \norm{\,  g \,  }_{H^k} ,
\end{align}
as soon as $\mathbb N \ni k > 1+d/2$ and analogously for all other appearing terms. Since 
$$
\norm{\,  g \,  }_{H^k} \leq K_3  \sum_{|j|+|\ell |\leq k}  \Big \| \, \partial_{x_\ell} \partial_{p_j} g \, \Big \|_\Lambda
$$ 
we obtain the desired estimate \eqref{qest}. 
\end{proof}

\medskip

{\bf Acknowledgement.} L.~Neumann thanks M. Escobedo for fruitful discussions on similar quantum kinetic models. 
C.~Sparber is thankful for the kind hospitality of the Johann Radon Institute for Applied Mathematics (RICAM).

\end{document}